\documentclass[aop,MSNbibl,seceqn,dvips]{arximspdf}
\usepackage{graphicx}
\doi{10.1214/12-AOP820} 
\volume{42}
\issue{4}
\pubyear{2014}
\firstpage{1644}
\lastpage{1665}

\makeatletter
\def\epsilon{\varepsilon}
\newcommand{\rrvert}{\vert}
\newcommand{\llvert}{\vert}
\newtheorem{theorem}{Theorem}[section]
\newtheorem{lemma}[theorem]{Lemma}
\newtheorem{proposition}[theorem]{Proposition}

\newproclaim{definition}[theorem]{Definition}
\newproclaim{remark}[theorem]{Remark}

\makeatother

\begin{document}
\begin{frontmatter}

\title{The Hausdorff dimension of the CLE gasket}
\runtitle{The Hausdorff dimension of the CLE gasket}

\begin{aug}
\author[A]{\fnms{Jason} \snm{Miller}\ead[label=u1,url]{http://jasonpmiller.org}},
\author[B]{\fnms{Nike} \snm{Sun}\corref{}\thanksref{t1}\ead[label=e2]{nikesun@stanford.edu}\ead[label=u2,url]{http://www-stat.stanford.edu/\textasciitilde nikesun/}}
\and
\author[C]{\fnms{David B.} \snm{Wilson}\ead[label=e3]{David.Wilson@microsoft.com}\ead[label=u3,url]{http://dbwilson.com}}
\thankstext{t1}{Supported in part by a Dept. of Defense (AFOSR) NDSEG Fellowship.}
\runauthor{J. Miller, N. Sun, and D. B. Wilson}
\affiliation{Massachussetts Institute of Technology, Microsoft
Research\break and Stanford University}
\address[A]{J. Miller\\
Department of Mathematics\\
Massachussetts Institute of Technology\\
77 Massachusetts Ave\\
Cambridge, Massachusetts 02139\\
USA}

\address[B]{N. Sun\\
Department of Statistics\\
Stanford University\\
Staford, California 94305\\
USA\\
\printead{u2}}

\address[C]{D. B. Wilson\\
Microsoft Research\\
One Microsoft Way\\
Redmond, Washington 98052\\
USA\\
\printead{u3}}
\end{aug}

\received{\smonth{6} \syear{2012}}
\revised{\smonth{10} \syear{2012}}

%
\begin{abstract}
The conformal loop ensemble $\mathrm{CLE}_\kappa$ is the canonical conformally
invariant probability measure on noncrossing loops in a proper simply
connected domain in the complex plane. The parameter $\kappa$ varies
between $8/3$ and $8$; $\mathrm{CLE}_{8/3}$ is empty while $\mathrm
{CLE}_8$ is a
single space-filling loop. In this work, we study the geometry of the
$\mathrm{CLE}$ \emph{gasket}, the set of points not surrounded by any
loop of
the $\mathrm{CLE}$. We show that the almost sure Hausdorff dimension
of the
gasket is bounded from below by $2-(8-\kappa)(3\kappa-8)/(32\kappa)$ when
$4<\kappa<8$. Together with the work of Schramm--Sheffield--Wilson
[\textit{Comm. Math. Phys.} \textbf{288} (2009) 43--53]
giving the upper bound for all $\kappa$ and the work of Nacu--Werner
[\textit{J. Lond. Math. Soc.} (2) \textbf{83} (2011) 789--809] giving the matching lower bound for $\kappa\le4$, this completes
the determination of the $\mathrm{CLE}_\kappa$ gasket dimension for
all values
of $\kappa$ for which it is defined. The dimension agrees with the
prediction of Duplantier--Saleur [\textit{Phys. Rev. Lett.} \textbf{63} (1989) 2536--2537] for the FK gasket.
\end{abstract}

%

\end{frontmatter}

\section{Introduction}
\label{sintro}

The conformal loop ensemble $\mathrm{CLE}_\kappa$ is the canonical
conformally invariant measure on countably infinite collections of
noncrossing loops in a proper simply connected domain $D$ in $\mathbb{C}$
\cite{MR2494457,SWCLE}. It is the loop analogue of $\mathrm
{SLE}_\kappa$, the
canonical conformally invariant measure on noncrossing paths. Whereas
$\mathrm{SLE}_\kappa$ arises as the scaling limit of a single macroscopic
interface of many two-dimensional discrete models \cite{MR1776084,MR2044671,MR2227824,MR2322705,MR2486487,arXiv10101356,MR2680496,chelkak-smirnov,CDCHKS}, $\mathrm{CLE}_\kappa$ describes the limit
of all of
the interfaces simultaneously. The parameter $\kappa$ varies between
$8/3$ and $8$; $\mathrm{CLE}_{8/3}$ is empty while $\mathrm{CLE}_8$
is a single
space-filling loop. $\mathrm{CLE}_\kappa$ for $\kappa\in(8/3,4]$
consists of
disjoint simple loops, while for $\kappa\in(4,8]$ the loops intersect
both themselves and each\vadjust{\vspace*{6pt}\eject} other (but are noncrossing). $\mathrm
{CLE}_3$ and
$\mathrm{CLE}_{16/3}$ are the scaling limits of the cluster boundaries
in the
square lattice critical Ising spin \cite{BDCHISINGCLE} and
FK-Ising~\cite{KempSmirFKISINGCLE} models, respectively, and
$\mathrm{CLE}
_6$ is the scaling limit of the cluster boundaries in critical
percolation on the triangular lattice \cite{MR2227824,MR2249794}.
$\mathrm{CLE}_4$ is the scaling limit of the level sets of the two-dimensional
discrete Gaussian free field~\cite{MSDGFFCLE}.\looseness=1

There are two different constructions of $\mathrm{CLE}_\kappa$. In
the first
construction, due to Werner \cite{MR2023758} and applicable for
$\kappa\in
[8/3,4]$, the loop ensemble is given by the outer boundaries of
Brownian loop soup clusters. In this paper, we make use of the second
construction, proposed by Sheffield \cite{MR2494457} and applicable for
$\kappa\in[8/3,8]$, based on branching $\mathrm{SLE}_\kappa(\kappa
-6)$. These
constructions have been proved equivalent for $\kappa\in[8/3,4]$
\cite
{SWCLE} (see also \cite{werner-wucle-exploration}).

\begin{figure}

\includegraphics{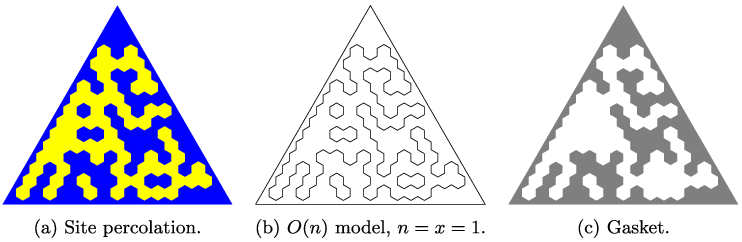}

\caption{Under the $O(n)$ model, a loop configuration $\omega$ has probability
proportional to $x^{e(\omega)} n^{\ell(\omega)}$
where $\ell(\omega)$ is the number of loops in $\omega$ and
$e(\omega)$ is the total length of all the loops. For $0\le n\le2$,
there is a critical value $x_c\equiv x_c(n)$ at which the $O(n)$ model
has a ``dilute phase,'' believed to converge to $\mathrm{CLE}_\kappa$ with
$n=-2\cos(4\pi/\kappa)$, $8/3\le\kappa\le4$. The $O(n)$ model at
$x>x_c$ is in a ``dense phase,'' again believed to converge to $\mathrm{CLE}
_\kappa$ with $n=-2\cos(4\pi/\kappa)$, but now with $4\le\kappa
\le
8$. Critical site percolation on the triangular lattice (left panel)
corresponds to the (dense phase) $O(n)$ model on the honeycomb lattice
with $n=x=1$ (center panel). Its gasket (right panel) is a
discretization of the $\mathrm{CLE}_6$ gasket.}
\label{fOn}
\end{figure}

%
%
%

\begin{figure}

\includegraphics{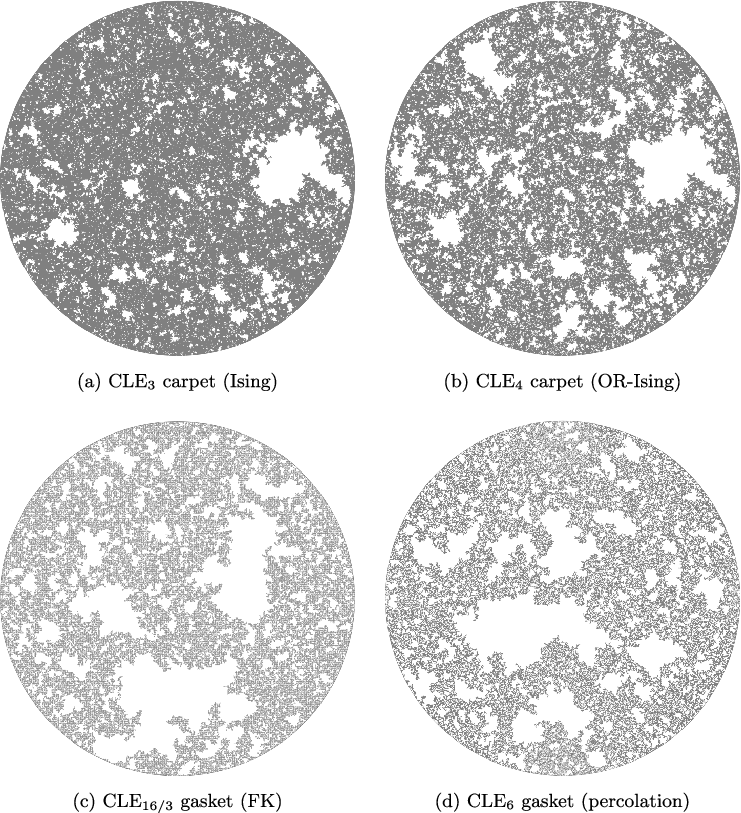}

\caption{Discrete simulations of the $\mathrm{CLE}_\kappa$ carpet
($\kappa\in
[8/3,4]$) or gasket
($\kappa\in(4,8]$) $\mathcal{G}_\kappa$ for $\kappa\in\{
3,4,16/3,6\}$. The
discretized $\mathcal{G}_\kappa$
(indicated in black above) is given by the set of points not surrounded
by any cluster boundary loop of a
discrete configuration sampled from a model known to converge to
$\mathrm{CLE}
_\kappa$. Note $\mathcal{G}_4\subseteq\mathcal{G}_3$
in our figures because the OR-Ising configuration used in \textup{(a)} is the binary OR of
two independent Ising configurations, one of which is used in \textup{(b)}.}
\label{fgasket}
\end{figure}

Let $\Gamma$ be a $\mathrm{CLE}_\kappa$ in $D$. The \emph{carpet}
($\kappa\in[8/3,4]$) or \emph{gasket} ($\kappa\in(4,8]$)
$\mathcal{G}$~of $\Gamma$ is the set of points not surrounded by any loop
of $\Gamma$. (In analogy with the Sierpi\'nski carpet and gasket, we
call $\mathcal{G}$ a carpet or gasket according to whether the loops
of $\Gamma$ are disjoint or intersecting, although occasionally we
loosely use gasket for both.) Since a.s. every neighborhood intersects
a loop, $\mathcal{G}$ is given equivalently by the closure of the
union of
the outermost loops of $\Gamma$. Figure~\ref{fOn} shows the gasket for
a discrete model, critical site percolation, that converges to $\mathrm{CLE}
_6$. Figure~\ref{fgasket} shows discrete simulations of $\mathcal{G}$
for $\kappa=3$ (Ising model), $\kappa=4$ (OR of two independent Ising
models, see \cite{SWCLE}, Proposition~10.2), $\kappa=16/3$ (FK-Ising
model), and $\kappa=6$ (critical percolation). The main result of this
article is the following theorem.\looseness=1\vadjust{\goodbreak}

%
\begin{theorem}
\label{thd}
Fix $\kappa\in(4,8)$ and let $\Gamma$ be a $\mathrm{CLE}_\kappa$
in a proper
simply connected domain $D$ in $\mathbb{C}$. Then with probability one the
Hausdorff dimension of the gasket $\mathcal{G}$ of $\Gamma$ is
%
%
\begin{equation}
\label{ehd} 2-\frac{(8-\kappa)(3\kappa-8)}{32\kappa}.
\end{equation}
\end{theorem}

The formula (\ref{ehd}) was first derived in the context of the
$O(n)$ model by Duplantier and Saleur \cite
{PhysRevLett632536,PhysRevLett64493}, who predicted the fractal dimension
of the $O(n)$ gasket (for $n\leq2$) using nonrigorous Coulomb gas
methods. The scaling limit of the $O(n)$ model is believed to be
$\mathrm{CLE}
_\kappa$, where $n=-2\cos(4\pi/\kappa)$ (\cite{MR2153402}, Conjecture~9.7, \cite{MR2494457}, Section~2.3). There are two
values of $\kappa$ associated to each $n<2$, corresponding to the
``dilute'' ($\kappa<4$) and ``dense'' ($\kappa>4$) phases of the $O(n)$
model. For further background see \cite{KN04}.

Schramm, Sheffield, and Wilson \cite{MR2491617} showed that for all
\mbox
{$8/3<\kappa<8$},\break (\ref{ehd})~gives the \emph{expectation
dimension}
of $\mathcal{G}$, the growth exponent of the expected number of balls of
radius $\epsilon$ needed to cover $\mathcal{G}$: this (a.s.) upper
bounds the
Minkowski dimension which in turn upper bounds the Hausdorff dimension.
(The expectation dimension for $\kappa=6$ was derived earlier by
Lawler, Schramm and Werner \cite{LSWone-arm}.) Nacu and Werner \cite
{MR2802511} used the Brownian loop soup construction to derive the
matching lower bound for the $\mathrm{CLE}_\kappa$ carpets ($\kappa
\le4$).

A lower bound on the Hausdorff dimension of a random fractal set is
obtained (by standard arguments) from a second moment estimate
controlling the probability that two given points lie near the set. The
complicated geometry of $\mathrm{CLE}$ loops prevents us from applying the
second moment method directly to $\mathcal{G}$, and instead we use a
``multi-scale refinement'' \cite{MR2123929}: we establish that with
arbitrarily small loss in the Hausdorff dimension we can restrict to
special classes of points in $\mathcal{G}$ whose correlation structure
\emph{at all scales} can be controlled.

\subsection*{Outline}

In Section~\ref{sprelim}, we review Sheffield's branching $\mathrm
{SLE}_\kappa
(\kappa-6)$ construction of $\mathrm{CLE}_\kappa$ [taking $\kappa
\in(4,8)$],
with an emphasis on its dependency structure. In Section~\ref{slbd},
we prove Theorem~\ref{thd}.

\section{Preliminaries}
\label{sprelim}

In this section, we review the exploration tree construction of
$\mathrm{CLE}
_\kappa$ for $\kappa\in(4,8)$ given in \cite{MR2494457} and then
collect several useful estimates for conformal maps.

\subsection{The continuum exploration tree}
\label{ssexplorationtree}

We begin by briefly recalling the definition of the $\mathrm
{SLE}_\kappa$ and
$\mathrm{SLE}_\kappa(\rho)$ processes. There are many excellent
surveys on the
subject (e.g., \cite{W03,MR2129588}) to which we refer the reader for
a more detailed introduction. The radial Loewner evolution in the unit
disk $\mathbb{D}$ is given by the differential equation
%
%
\begin{equation}
\label{elderadial} \dot g_t(z) = -g_t(z)
\frac{g_t(z) + W_t}{g_t(z) - W_t},\qquad g_0(z)=z,
\end{equation}
where $W_t$ is a continuous function which takes values in $\partial
\mathbb{D}$.
We refer to $W_t$ as the \emph{driving function} of the Loewner
evolution. For $z\in\mathbb{D}$, let
\[
T^z\equiv\sup\bigl\{t \geq0 \dvtx \bigl|g_t(z)\bigr| < 1\bigr\}
\]
and
\[
K_t\equiv\bigl\{z\in\mathbb{D}\dvtx T^z\le t\bigr\}.
\]
For each $t\ge0$, $g_t$ is the unique conformal transformation
$\mathbb{D}
\setminus K_t\to\mathbb{D}$ with $g_t(0)=0$ and $g_t'(0)>0$. The (random)
growth process $(K_t)_{t\ge0}$ associated with $W_t=\exp(i\sqrt {\kappa
}B_t)$, where $B_t$ is a standard Brownian motion, is the radial
$\mathrm{SLE}
_\kappa$ process introduced by Schramm \cite{MR1776084}. Time is
parametrized by negative log-conformal radius, that is, $g_t'(0)=e^t$.
It was proved by Rohde and Schramm \cite{MR2153402} ($\kappa\ne8$)
and Lawler, Schramm, and Werner \cite{MR2044671} (\mbox{$\kappa= 8$})
that there is a \emph{curve} $\eta\dvtx [0,\infty)\to\overline\mathbb
{D}$ starting
at $\eta(0)=1$ such that $\mathbb{D}\setminus K_t$ is the unique connected
component of $\mathbb{D}\setminus\eta[0,t]$ containing $0$: we say that
$\eta$ \emph{generates} the process $K_t$ and call $\eta$ the
radial $\mathrm{SLE}_\kappa$ \emph{trace}. In this setting,
$W_t=\lim_{z\to
\eta(t)} g_t(z)$, where the limit is taken with $z\in\mathbb
{D}\setminus
K_t$. For $\kappa<8$, Lawler \cite{lawlerendpoint} proved that $\lim_{t\to\infty} \eta(t)=0$, so $\eta\dvtx [0,\infty]\to\overline\mathbb
{D}$ defines
a curve traveling from $\eta(0)=1$ to $\eta(\infty)=0$ in $\overline
\mathbb{D}$.

Let $D$ be a proper simply connected domain in $\mathbb{C}$. For any conformal
transformation $f\dvtx \mathbb{D}\to D$, we take the image of radial
$\mathrm{SLE}_\kappa$
in $\mathbb{D}$ under $f$ to be the definition of radial $\mathrm
{SLE}_\kappa$ in $D$
from $f(1)$ to $f(0)$, with $f(1)$ interpreted as a prime end. If $f$
extends continuously to $\overline\mathbb{D}$ (equivalently if
$\partial D$ is given by
a closed curve, see \cite{MR1217706}, Theorem~2.1), then radial $\mathrm{SLE}
_\kappa$ in $D$ is a.s. a continuous curve. It was proved by Garban,
Rohde and Schramm \cite{MR2891697} that radial $\mathrm{SLE}_\kappa$
with $\kappa<8$ in a \emph{general} proper simply connected domain is
a.s. continuous except possibly at its starting point.

We now describe the radial $\mathrm{SLE}_\kappa(\rho)$ processes, a natural
generalization of radial $\mathrm{SLE}_\kappa$ first introduced in
\cite{MR1992830}, Section~8.3. For $w,o\in\partial\mathbb{D}$, radial
$\mathrm{SLE}_\kappa
(\rho)$ with starting configuration $(w,o)$ is the (random) growth
process associated with the solution of (\ref{elderadial}) where the
driving function solves the SDE
%
%
\begin{equation}
\label{eslekr} dW_t = -\frac{\kappa}{2} W_t \,d t + i
\sqrt{\kappa} W_t \,d B_t - \frac{\rho}{2}
W_t \frac{W_t+O_t}{W_t-O_t} \,d t,\qquad W_0=w
\end{equation}
with $O_t=g_t(o)$, the \emph{force point}. It is easy to see that
(\ref{eslekr}) has a unique solution up to time $\tau_=\equiv\inf
\{t\ge0\dvtx W_t=O_t\}$.

The weight $\rho=\kappa-6$ is special because it arises as a
coordinate change of ordinary chordal $\mathrm{SLE}_\kappa$ from $w$
targeted at
$o$. A consequence is that \emph{radial $\mathrm{SLE}_\kappa(\kappa
-6)$ is target
invariant:} radial $\mathrm{SLE}_\kappa(\kappa-6)$ in $\mathbb
{D}$ with starting
configuration $(w,o)$ and target $a\in\mathbb{D}$ has the same law (modulo
time change) as an ordinary chordal $\mathrm{SLE}_\kappa$ in $\mathbb
{D}$ from $w$ to
$o$, up to the first time the curve disconnects $a$ and $o$ \cite
{MR2188260}.

We now explain how to construct a solution to (\ref{eslekr}) which
is defined even after time $\tau_=$. A more detailed treatment is
provided in \cite{MR2494457}, Section~3; we give here a brief summary
following \cite{MR2491617}. For $\rho>-\kappa/2-2$, there is a random
continuous process $\theta_t$ taking values in $[0,2\pi]$ which
evolves according to the SDE
%
%
\begin{equation}
\label{etheta} d\theta_t =\sqrt{\kappa} \,d B_t+
\frac{\rho+2}{2}\cot(\theta_t/2) \,d t
\end{equation}
on each interval of time for which $\theta_t\notin\{0,2\pi\}$, and
is instantaneously reflecting at the endpoints, that is, the set $\{
t\dvtx \theta_t\in\{0,2\pi\}\}$ has Lebesgue measure zero. (This
diffusion was studied in \cite{LSWone-arm} for $\rho=0$.) In other
words, $\theta_t$ is a random continuous process adapted to the
filtration of $B_t$ which a.s. satisfies
\[
\partial_t[\theta_t-\sqrt{\kappa}B_t]=
\frac{\rho+2}{2}\cot (\theta_t/2)
\]
for all $t$ for which the right-hand side is finite. The law of this
process is uniquely determined by $\theta_0$, and moreover the process is
pathwise unique \cite{MR2494457}, Proposition~4.2. It then follows from the
strong Markov property of Brownian motion that $\theta_t$ has the
strong Markov property.

When $\rho\ge\kappa/2-2$, the $\theta_t$ process governed by
SDE (\ref{etheta}) is repelled so strongly by $0$ and $2\pi$ that it
almost surely never reaches either endpoint. When $\rho=-2$ the
diffusion $\theta_t$ is simply reflected Brownian motion. When $\rho
<-2$, the $\theta_t$ process is attracted to the singularity and its
analysis requires more care, but it still makes sense when $\rho
>-\kappa/2-2$ \cite{MR2494457,MR2491617}. When $\rho\leq-\kappa
/2-2$, the
$\theta_t$ process is attracted so strongly to the endpoints that once
it hits either one it remains glued there. In the intermediate regime,
$-\kappa/2-2<\rho<\kappa/2-2$, the $\theta_t$ process hits the
endpoints $0$ and $2\pi$, but is instantaneously reflecting. When
$\rho=\kappa-6$, this corresponds to the range $8/3<\kappa<8$.

We then set
%
%
\begin{equation}
\label{eslekrint} \arg W_t =\arg w+\sqrt{\kappa}B_t+
\frac{\rho}{2}\int_0^t\cot(
\theta_s/2) \,d s.
\end{equation}
That the above integral is a.s. finite follows by the comparison of
$\theta_t/\sqrt{\kappa}$ [resp., $(2\pi-\theta_t)/\sqrt{\kappa}$]
with a $\delta$-dimensional Bessel process, as described above; see, for
example, the proof of Lemma~\ref{lcond}. We then \emph{define}
radial $\mathrm{SLE}_\kappa(\rho)$ in $\mathbb{D}$ with starting
configuration $(w,o)$
to be the solution to (\ref{elderadial}) with driving function $W_t$
defined by~(\ref{eslekrint}). The force point $O_t\equiv g_t(o)$
satisfies $W_t=O_t e^{i\theta_t}$, and we interpret $\theta_t=0$ as
$O_t=W_t e^{i0^-}$ ($\arg O_t$ just below $\arg W_t$) and similarly
$\theta_t=2\pi$ as $O_t=W_t e^{i 0^+}$. For $\rho\ge\kappa/2-2$, the
laws of radial $\mathrm{SLE}_\kappa(\rho)$ and ordinary radial
$\mathrm{SLE}_\kappa$ are
mutually absolutely continuous up to any fixed positive time, so
$\mathrm{SLE}
_\kappa(\rho)$ is a.s. generated by a curve by the result of \cite
{MR2153402}. In \cite{MSIMAG}, it is established that $\mathrm
{SLE}_\kappa
(\rho)$ is a.s. generated by a curve for all $\rho>-2$ (see
Remark~\ref{rcle}); when $\rho=\kappa-6$ this corresponds to $\kappa
>4$. Radial $\mathrm{SLE}_\kappa(\rho)$ in a general proper simply connected
domain is defined again by conformal transformation, but the analogue
of the continuity result of \cite{MR2891697} is not known
for $\rho\ne0$.

The target invariance of radial $\mathrm{SLE}_\kappa(\kappa-6)$ processes
continues to hold after time $\tau_=$, and from this we can construct a
coupling of radial $\mathrm{SLE}_\kappa(\kappa-6)$ processes
targeted at a
countable dense subset of $\mathbb{D}$.

\begin{proposition}[({\cite{MR2494457}, Proposition~3.14 and Section~4.2})]
\label{ptargetinv}
Let $(a_k)_{k\in\mathbb{N}}$ be a countable dense sequence in
$\mathbb{D}$. For
$4<\kappa<8$, there exists a coupling of radial $\mathrm{SLE}_\kappa
(\kappa-6)$
curves $\eta^{a_k}$ in $\mathbb{D}$ from $1$ to $a_k$ started from
$(w,o)=(1,1e^{i0^-})$ such that for any $k,\ell\in\mathbb{N}$, $\eta^{a_k}$
and $\eta^{a_\ell}$ agree a.s. (modulo time change) up to the first
time that the curves separate $a_k$ and $a_\ell$ and evolve
independently thereafter.
\end{proposition}
(For $8/3<\kappa<4$, the $\mathrm{SLE}_\kappa(\kappa-6)$ traces are
not known to
be curves, which makes the corresponding statement in this case more
complicated. The case $\kappa=4$ is special, and was dealt with
separately by Sheffield \cite{MR2494457}.)

From the coupling $(\eta^{a_k})_{k \in\mathbb{N}}$ defined in
Proposition~\ref{ptargetinv}, we can a.s. uniquely define (modulo
time change) for each $a\in\overline\mathbb{D}$ a curve $\eta^a$
targeted at $a$,
by considering a subsequence $(a_{k_n})$ converging to $a$. Then $\eta
^a$ is a radial \mbox{$\mathrm{SLE}_\kappa(\kappa-6)$}, and we write
$\theta
^a_t,W^a_t,O^a_t$ for the corresponding processes of (\ref{etheta})
and (\ref{eslekrint}). The complete collection of curves $(\eta
^a)_{a\in\overline\mathbb{D}}$ is the \emph{branching $\mathrm
{SLE}_\kappa(\kappa-6)$} or
\emph{continuum exploration tree} of \cite{MR2494457}.

\subsection{Loops from exploration trees}\label{sstreesloops}

For $4<\kappa<8$,
the $\mathrm{CLE}_\kappa$ loops $\mathcal{L}^a$ surrounding $a\in
\mathbb{D}$ are defined in
terms of the branch $\eta^a$ of the exploration tree as follows:
\begin{longlist}[1.]
\item[1.] Let $\tau_{\mathrm{ccw}}^a\equiv\inf\{t\ge0\dvtx \theta^a_t=2\pi
\}$, the
first time $\eta^a$ forms a counterclockwise loop surrounding $a$.
\item[2.] If $\tau_{\mathrm{ccw}}^a=\infty$, then there are no loops
surrounding $a$
and we set $\mathcal{L}^a$ to be the empty sequence. If $\tau
_{\mathrm{ccw}}^a<\infty$
let $\acute\tau_{\mathrm{ccw}}^a\equiv\sup\{t<\tau_{\mathrm
{ccw}}^a\dvtx \theta^a_t=0\}$, let
$o^a\equiv\eta^a(\acute\tau_{\mathrm{ccw}}^a)$, and let
$\widetilde\eta^a$ be the branch
$\eta^{o^a}$, reparametrized so that $\widetilde\eta^a|_{[0,\tau
_{\mathrm{ccw}}^a]}=\eta
^a|_{[0,\tau_{\mathrm{ccw}}^a]}$. The outermost loop $\mathcal
{L}^a_1$ surrounding $a$ is
defined to be $\widetilde\eta^a|_{[\acute\tau_{\mathrm
{ccw}}^a,\infty]}$.
\end{longlist}

\begin{figure}

\includegraphics{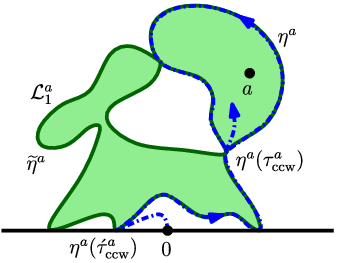}

\caption{Branching $\mathrm{SLE}_\kappa(\kappa-6)$ construction of $\mathrm
{CLE}_\kappa$ ($4<\kappa
<8$) process $\Gamma$ in $\mathbb{H}$.
For each $a\in\mathbb{H}$, $\eta^a$ (dashed blue line) is the branch
of the
exploration tree targeted at $a$.
It evolves as a radial $\mathrm{SLE}_\kappa(\kappa-6)$ which,
whenever it hits the
domain boundary or its past hull,
continues in the complementary connected component containing $a$. Let
$\tau_{\mathrm{ccw}}^a$ be the first time
$t$ that $\eta^a$ completes a counterclockwise loop surrounding $a$;
the location of the force point
at time $\tau_{\mathrm{ccw}}^a$ is $o^a\equiv\eta^a(\acute\tau
_{\mathrm{ccw}}^a)$ for some
$\acute\tau_{\mathrm{ccw}}^a<\tau_{\mathrm{ccw}}^a$. The outermost
loop $\mathcal{L}_1^a$ of $\Gamma$ containing $a$ is $\eta
^{o^a}|_{[\acute\tau_{\mathrm{ccw}}^a,\infty]}$. Successive
loops are defined in analogous fashion. $\mathcal{L}^a_1$ is necessarily
counterclockwise and pinned
at $\eta^a(\acute\tau_{\mathrm{ccw}}^a)$. It is disjoint from the
domain boundary if
and only if $a$ is first surrounded by a clockwise loop.}
\label{fcletree}
\end{figure}

\begin{figure}[b]

\includegraphics{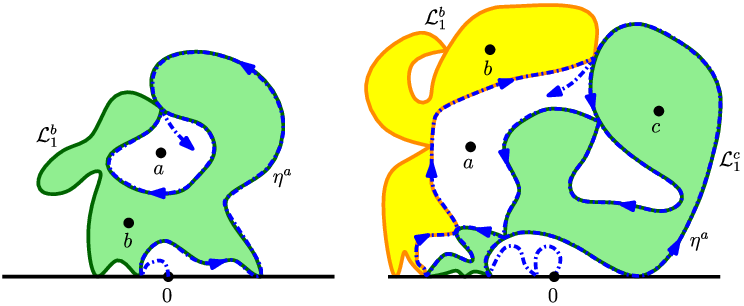}

\caption{Clockwise loops of $\eta^a$ (dashed blue line) are not
$\mathrm{CLE}
$ loops, but correspond
either to complementary connected components of $\mathrm{CLE}$ loops (left
panel) or complementary
connected components of chains of $\mathrm{CLE}$ loops (right panel). The
$\mathrm{CLE}$ process is renewed within each clockwise loop
(Proposition \protect\ref{pcle}).}
\label{fcletree2}
\end{figure}

If $\mathcal{L}^a_1$ is defined, it is necessarily counterclockwise and
pinned at $\eta^a(\acute\tau_{\mathrm{ccw}}^a)$, and for any point
$b$ surrounded by
$\mathcal{L}^a_1$ we have $\mathcal{L}^b_1=\mathcal{L}^a_1$.
Moreover, $\eta^a(\acute\tau_{\mathrm{ccw}}^a)$ lies on $\partial
\mathbb{D}$ if and only if
$\eta^a$ has not previously made a clockwise loop around $a$
\cite{MR2494457}, Lemma~5.2. The next loop $\mathcal{L}^a_2$ surrounding $a$
is then
defined in analogous fashion, and continuing in this way gives the full
$\mathrm{CLE}_\kappa$ process $\Gamma$ in $\mathbb{D}$. See
Figures~\ref{fcletree}
and~\ref{fcletree2}.

\begin{remark}
\label{rcle}
For $4<\kappa<8$, assuming the conjecture that chordal\break  $\mathrm
{SLE}_\kappa(\kappa
-6)$ processes are generated by continuous curves with reversible law
\cite{MR2494457}, Conjecture 3.11, it was shown \cite{MR2494457}, Proposition~5.1 and
Theorem~5.4, that $\mathrm{CLE}_\kappa$ loops are continuous,
and that the
law of the full ensemble is independent of the choice of root for the
exploration tree. This conjecture was proved in works of Miller and
Sheffield (\cite{MSIMAG}, Theorem~1.3 and \cite{MSIMAG3}, Theorems 1.1
and 1.2), so these properties hold. (The analogous
continuity and root-invariance statements are immediate for $\kappa\in
[8/3,4]$ by the equivalence of $\mathrm{CLE}_\kappa$ and the outer
boundaries of
loop soups~\cite{SWCLE}; see also \cite{werner-wucle-exploration}.)
\end{remark}

The $\mathrm{CLE}_\kappa$ process in a general proper simply
connected domain is
defined by conformal transformation, so the law of $\mathrm
{CLE}_\kappa$ is
conformally invariant. Moreover, conditional on the collection of all
of the outermost loops, the law of the loops contained in the connected
component $D^a$ of $\mathbb{D}\setminus\mathcal{L}_1^a$
containing $a$ is equal to
that of a $\mathrm{CLE}_\kappa$ in $D^a$ independently of the loops
of $\Gamma$
which are not contained in $D^a$. The key observation which we use to
prove Theorem~\ref{thd} is that there are additional sources of
conditional independence in $\mathrm{CLE}_\kappa$ when $\kappa>4$,
in particular:

\begin{proposition}\label{pcle}
Suppose $z\in D$ is surrounded by a clockwise loop $\mathcal{C}$ in the
$\mathrm{SLE}_\kappa(\kappa-6)$ exploration tree of $D$ (as in
Figure~\ref{fcletree2}), allowing the domain boundary to form part of the loop
$\mathcal{C}$. If $U$ is the connected component of $D\setminus
\mathcal{C}$
containing~$z$, then the law of the $\mathrm{CLE}_\kappa$ loops contained
within $\overline U$ is that of a $\mathrm{CLE}_\kappa$ in $U$, independent
of the $\mathrm{CLE}_\kappa$ loops outside of $U$.
\end{proposition}

The $\mathrm{SLE}_\kappa(\kappa-6)$ exploration tree for $\kappa>4$
has such
clockwise loops, which are not $\mathrm{CLE}$ loops, and so provide additional
renewal events.\

\subsection{Diffusion estimate}

\begin{proposition}[({\cite{MR2491617}, equation (4)})]\label{pssw}
Suppose $8/3<\kappa<8$, and let $\theta_t$ be the process defined
above started from $\theta_0=0$, evolving according to SDE (\ref
{etheta}) in $(0,2\pi)$ and instantaneously reflecting at the
endpoints $\{0,2\pi\}$. Then $\mathbb{P}[\theta_s< 2\pi\ \forall
s\le t]
\asymp e^{-\alpha t}$ where
%
%
\begin{equation}
\label{alpha} \alpha\equiv\frac{(8-\kappa)(3\kappa-8)}{32\kappa}.
\end{equation}
\end{proposition}

It is this diffusion exponent $\alpha$ which gave rise to the result
of \cite{MR2491617} that the gasket has expectation dimension
$2-\alpha$,
implying an upper bound of $2-\alpha$ for the Hausdorff dimension, for
which Theorem~\ref{thd} provides the matching lower bound. The actual
value of $\alpha$ does not play a significant role in the proof of
Theorem~\ref{thd}, except that we use $0<\alpha<2$. (Of course,
$\alpha\leq2$ is a necessary condition for showing that the Hausdorff
dimension is $2-\alpha$.)

\subsection{Distortion estimates}
\label{ssdistort}

For a proper simply connected domain $D$ and $w\in D$, let $\operatorname{CR}
(w,D)$ denote the conformal radius of $D$ with respect to $w$, that is,
$\operatorname{CR}(w,D)\equiv f'(0)$ for $f$ the unique conformal map
$\mathbb{D}\to
D$ with $f(0)=w$ and $f'(0)>0$. Let $\operatorname{rad}(w,D)\equiv\inf\{
r\dvtx B_{r}(w)\supseteq D\}$ denote the out-radius of $D$ with respect
to $w$. By the Schwarz lemma and the Koebe one-quarter theorem,
%
%
\begin{equation}
\label{eircror} \operatorname{dist}(w,\partial D) \le\operatorname{CR}(w,D) \le\bigl[4
\operatorname{dist}(w,\partial D)\bigr] \wedge\operatorname{rad}(w,D).
\end{equation}
Further (see, e.g., \cite{MR1217706}, Theorem~1.3)
\[
\frac{|\zeta|}{(1+|\zeta|)^2} \leq\frac{|f(\zeta)-w|}{\operatorname{CR}(w,D)} \leq\frac{|\zeta|}{(1-|\zeta|)^2}.
\]
As a consequence,
%
%
\begin{equation}
\label{egrowth} \frac{|\zeta|}{4} \le\frac{|f(\zeta)-w|}{\operatorname{CR}(w,D)} \le4|\zeta|,
\end{equation}
where the right-hand inequality holds for $|\zeta|\le1/2$.

\section{Proofs}
\label{slbd}

Recall that a $\mathrm{CLE}_\kappa$ process in a general simply
connected domain
$D$ is defined as the image under a conformal transformation $f\dvtx \mathbb
{D}\to
D$ of a $\mathrm{CLE}_\kappa$ process $\Gamma$ in $\mathbb{D}$.
Since $f|_{r \mathbb{D}}$ for
any $0<r<1$ is bi-Lipschitz and so preserves Hausdorff dimension, and
the Hausdorff dimension of a countable union is the supremum of the
Hausdorff dimensions, we see that $f$ preserves Hausdorff dimension,
and so it suffices to prove Theorem~\ref{thd} with $D=\mathbb{D}$.
Thus, for
the remainder $\Gamma$ denotes a $\mathrm{CLE}_\kappa$ ($4<\kappa
<8$) process on
$\mathbb{D}$, constructed from the collection of radial $\mathrm
{SLE}_\kappa(\kappa-6)$
curves $(\eta^z)_{z\in\mathbb{D}}$ jointly defined on a probability space
$(\Omega,\mathcal{F},\mathbb{P})$, as given by the remark following
Proposition~\ref
{ptargetinv}. In Section~\ref{ssevents}, we define our multi-scale
refinement of the gasket $\mathcal{G}$ of $\Gamma$, and state the main
result of the section, the second moment estimate Lemma~\ref{lind} on
the correlation structure of the set of ``perfect points'' identified
by the refinement. We then use the $\mathrm{CLE}$ renewal property of
Proposition~\ref{pcle} to reduce Lemma~\ref{lind} to a lower bound
on the probability of a single event. This bound is given by
Proposition~\ref{pexp}, which we prove in Section~\ref{sscw}. The
Hausdorff dimension lower bound follows from Lemma~\ref{lind} by
standard arguments which we give in Section~\ref{sshd}, thereby
concluding the proof of Theorem~\ref{thd}.

\begin{figure}

\includegraphics{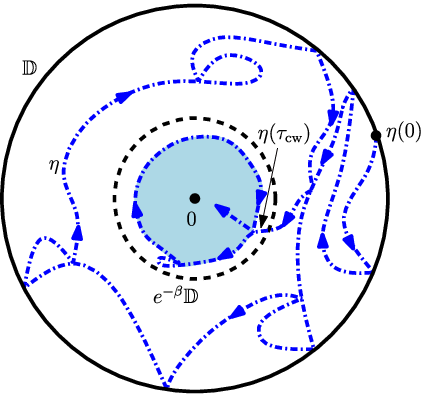}

\caption{A single level of the multi-scale argument we use to prove
the lower bound of Theorem~\protect\ref{thd}.
The curve $\eta$ (dashed blue line) is a radial $\mathrm{SLE}_\kappa
(\kappa-6)$
targeted at zero. For $\beta>0$, $E(\eta)$ is the
event that the first time $\tau_{\mathrm{cw}}$ that $\eta$ closes a
clockwise loop
$\mathcal{C}$ surrounding $0$ with $\mathcal{C}\subset e^{-\beta}
\mathbb{D}$
is finite, and further that $\eta$ makes no counterclockwise loop
surrounding $0$ before $\tau_{\mathrm{cw}}$. On the event $E(\eta)$,
set $D(\eta)$ (light blue region) to be the connected component of
$\mathbb{D}
\setminus\eta[0,\tau_{\mathrm{cw}}]$ containing $0$.}
\label{fevt}
\end{figure}

\subsection{Clockwise loops in small disks}
\label{ssevents}

We now describe our multi-scale refinement of the gasket $\mathcal{G}$
which identifies a subset of ``perfect points'' (following the
terminology of \cite{MR2123929,MR2642894}) in $\mathcal{G}$,
satisfying a
certain restriction at all scales which makes their correlation
structure easy to analyze. Let $\beta>0$ be a parameter (which we will
send to $\infty$). Let $\eta$ be any curve defined on $(\Omega
,\mathcal{F},\mathbb{P}
)$ and traveling in $\overline\mathbb{D}$ from $\partial\mathbb{D}$
to $0$. For any
curve $\eta$ defined on $(\Omega,\mathcal{F},\mathbb{P})$ and
traveling in $\overline\mathbb{D}$
from $\partial\mathbb{D}$ to $0$, define $E(\eta)\subseteq\Omega$
to be the event that
\begin{longlist}[(ii)]
\item[(i)] the first time $\tau_{\mathrm{cw}}$ that $\eta$ closes a
clockwise loop
$\mathcal{C}$ surrounding $0$ with $\mathcal{C}\subset e^{-\beta}
\mathbb{D}$ is
finite; and
\item[(ii)]$\eta$ makes no counterclockwise loop surrounding $0$ before
time $\tau_{\mathrm{cw}}$.
\end{longlist}
On the event $E(\eta)$, set $D(\eta)$ to be the connected component
of $\mathbb{D}\setminus\mathcal{C}$ containing the origin. See
Figure~\ref{fevt} for an illustration.

We then define events $E_j$ and domains $D_j\ni0$, both nonincreasing
in $j$ for $j\ge0$, as follows: let $(E_0,D_0)\equiv(\Omega,\mathbb
{D})$, and
suppose inductively that $(E_j,D_j)$ has been defined. Let $\mathsf{g}_j$ be
the uniformizing map $D_j\to\mathbb{D}$ with $\mathsf{g}_j(0)=0$ and $\mathsf{g}
_j'(0)>0$, and let
\[
\tau_j\equiv\inf\bigl\{t \geq0\dvtx \eta(t)\in D_j\bigr
\},\qquad \mathsf{g}_j\eta\equiv\bigl(\mathsf{g}_j
\eta(\tau _j+s)\bigr)_{s\ge0}.
\]
We then set
\[
E_{j+1}(\eta) \equiv E_j(\eta) \cap E(\mathsf{g}_j\eta ),\qquad D_{j+1}(\eta) \equiv\mathsf{g}_j^{-1}D(\mathsf{g}_j\eta).
\]
For $z\in\mathbb{D}$, let
\[
\psi(\zeta)\equiv\psi_z(\zeta)\equiv\frac{\zeta-z}{1-\bar
z\zeta},
\]
the conformal automorphism of $\mathbb{D}$ with $\psi(z)=0$ and \mbox
{$\psi
'(z)=(1-|z|^2)^{-1}>0$}. For $\eta^z$, the branch of the $\mathrm{SLE}
_\kappa(\kappa-6)$ exploration tree targeted at $z$, we set
%
%
\begin{equation}
\label{eDzn} E^z_j \equiv E_j\bigl(
\psi_z\eta^z\bigr),\qquad D^z_j\equiv
\psi_z^{-1}D_j\bigl(\psi_z
\eta^z\bigr).
\end{equation}
The \emph{perfect points} in the multi-scale refinement of the
gasket $\mathcal{G}$ are the points $z\in\mathbb{D}$ for which
$\bigcap_{j\ge
0}E^z_j$ occurs. The main estimate needed to lower bound the Hausdorff
dimension is the following estimate on their correlation structure.
%
\begin{lemma}
\label{lind}
For sufficiently large $\beta$ there exists $\varepsilon\equiv
\varepsilon(\beta
)<\infty$ with $\lim_{\beta\to\infty}\varepsilon(\beta)=0$ such
that for
all $z,w\in\mathbb{D}$,
\[
\frac{\mathbb{P}[E^z_n \cap E^w_n]}{\mathbb{P}[E^z_n]\mathbb{P}[E^w_n]} \le \biggl(\frac{e^{\beta}}{|z-w|} \biggr)^{\alpha(1+\varepsilon)},
\]
where $\alpha$ is given by (\ref{alpha}).
\end{lemma}

In the remainder of this subsection, we reduce the proof of this lemma
to a lower bound on the probability of the event $E^0_1$,
Proposition~\ref{pexp}, which we prove in Section~\ref{sscw}. We
begin with some easy estimates comparing the domains $D^z_j$ to disks
$B_{e^{-j\beta}}(z)$.

\begin{lemma}\label{ldistort}
For $\beta\ge\log2$, $j\geq1$, and $z\in\mathbb{D}$,
%
%
\begin{equation}
\label{edistortionbound} \operatorname{rad}\bigl(z,D^z_j\bigr)
\le8 e^{-j\beta} \qquad\mbox{on } E^z_j.
\end{equation}
\end{lemma}

\begin{pf}
We first consider the domains $D_j\equiv D_j(\eta)$, defined on the
event $E_j(\eta)$, for any curve $\eta$ traveling in $\overline
\mathbb{D}$ from
$\partial\mathbb{D}$ to $0$. [We will later take $\eta= \psi_z(\eta
^z)$, where
$\eta^z$ is the branch of the $\mathrm{SLE}_\kappa(\kappa-6)$ exploration
tree targeted at $z$, and $\psi_z$ is the M\"obius transformation
defined above which maps $z$ to $0$.] Recall the definition of the
uniformizing map $\mathsf{g}_j\dvtx D_j\to\mathbb{D}$. By the
definition of $D_j$ and
by (\ref{eircror}),
%
%
\begin{equation}
\label{econfradrad} \operatorname{CR}(0,\mathsf{g}_{j-1}D_j)
\le\operatorname{rad}(0,\mathsf{g}_{j-1}D_j) \le
e^{-\beta}.
\end{equation}
Since
\[
\operatorname{CR}(0,\mathsf{g}_{j-1}D_j)=
\frac{1}{(\mathsf{g}_j\circ\mathsf{g}
_{j-1}^{-1})'(0)} = \frac{\operatorname{CR}(0,D_j)}{\operatorname{CR}(0,D_{j-1})},
\]
we have that
%
%
\begin{equation}
\label{eprodcr} \operatorname{CR}(0,D_j) =\prod
_{\ell=1}^j\operatorname{CR}(0,g_{\ell-1}D_\ell)
\le e^{-j\beta}.
\end{equation}
Applying the right-hand inequality of (\ref{egrowth}) with $f=\mathsf{g}_{j-1}^{-1}$ gives
\[
\frac{|\zeta|}{\operatorname{CR}(0,D_{j-1})} \le4\bigl|\mathsf{g}_{j-1}(\zeta)\bigr| \le4
e^{-\beta}\qquad \mbox{when }\zeta\in\partial D_j,
\]
using that $\zeta\in\partial D_j$ implies $|\mathsf{g}_{j-1}(\zeta)|\le
e^{-\beta}\le1/2$. Rearranging and combining with (\ref{eprodcr}) gives
%
%
\begin{equation}
\label{ezerodistort} |\zeta|\leq4 e^{-\beta} \operatorname{CR}(0,D_{j-1})
\le4 e^{-j\beta} \qquad\mbox{when }\zeta\in\partial D_j.
\end{equation}

For any $z\in\mathbb{D}$, (\ref{ezerodistort}) is satisfied with
$D_j=D_j(\psi_z\eta^z)=\psi_z D^z_j$ on the event $E^z_j$.
We have $\psi_z^{-1}(\zeta)=(z+\zeta)/(1+\bar z\zeta)$, so
\[
\bigl\llvert \psi_z^{-1}(\zeta)-z\bigr\rrvert =\biggl
\llvert \zeta\frac{1-|z|^2}{1+\bar z\zeta} \biggr\rrvert \le|\zeta| \frac{1-|z|^2}{1-|z|} = |
\zeta|\bigl(1+|z|\bigr) \le2|\zeta|,
\]
giving $\operatorname{rad}(z,D^z_j) \le2 \operatorname{rad}(0,\psi_z D^z_j)
\leq8
e^{-j\beta}$ as claimed.
\end{pf}

Let $\mathcal{F}^z_j$ denote the $\sigma$-algebra generated by $\eta
^z$ up to
the time $\tau^z_j$ that $\eta^z$ closes the clockwise loop forming
the boundary of $D^z_j$ (if $E^z_j$ does not occur then $\tau^z_j =
\infty$). By the conformal Markov property of radial $\mathrm
{SLE}_\kappa
(\rho)$, for $m\le n$, we have
\[
\mathbb{P}\bigl[E^z_n |\mathcal{F}_m^z
\bigr] \mathbf{1}_{E^z_m} = \mathbb{P}\bigl[E^z_{n-m}
\bigr] \mathbf{1}_{E^z_m} = \mathbb{P}\bigl[E^0_{n-m}
\bigr] \mathbf{1}_{E^z_m},
\]
and consequently
\[
\mathbb{P}\bigl[E^z_n\bigr] =\mathbb{E}\bigl[\mathbb{P}
\bigl[E^z_n |\mathcal{F}_{n-1}^z
\bigr]\mathbf{1}_{E^z_{n-1}}\bigr] =\mathbb{P}\bigl[E^0_1
\bigr]\mathbb{P}\bigl[E^z_{n-1}\bigr] =\cdots=\mathbb{P}
\bigl[E^0_1\bigr]^n.
\]

\begin{proposition}\label{pexp}
$\!\!\!$There exists a constant $c > 0$ such that $\mathbb{P}[E^0_1]\ge(c
e^{\alpha\beta
})^{-1}$ for sufficiently large $\beta$, where $\alpha$ is given by
(\ref{alpha}).
\end{proposition}

The proof of this proposition is deferred to Section~\ref{sscw}, but
we show now how to use it to deduce Lemma~\ref{lind}.

\begin{pf*}{Proof of Lemma~\ref{lind}}
Given $z,w\in\mathbb{D}$, let $m\in\mathbb{N}$ be defined by $8
e^{-m\beta}<|z-w|\le8  e^{-(m-1)\beta}$. If $E^z_m\cap E^w_m$
occurs, then Lemma~\ref{ldistort} implies $w \notin D_m^z$ which in
turn implies $D_m^z \cap D_m^w = \varnothing$. So for $n\geq m$, $E^z_n$
and $E^w_n$ are conditionally independent given $E^z_m \cap E^w_m$, and
in fact
\[
\mathbb{P}\bigl[E^z_n \cap E^w_n
|E^z_m \cap E^w_m\bigr] =
\mathbb{P}\bigl[E^0_{n-m}\bigr]^2.
\]
Therefore
\begin{eqnarray*}
\mathbb{P}\bigl[E^z_n \cap E^w_n
\bigr] &=& \mathbb{P}\bigl[E^z_m \cap E^w_m
\bigr] \mathbb{P}\bigl[E^0_{n-m}\bigr]^2
\\
&\le&\frac{(\mathbb{P}[E^0_m]\mathbb{P}[E^0_{n-m}])^2}{\mathbb{P}[E^0_m]} = \frac{\mathbb{P}[E^0_n]^2}{\mathbb{P}[E^0_1]^m} \le\bigl(ce^{\alpha\beta}
\bigr)^m \mathbb{P}\bigl[E^0_n
\bigr]^2,
\end{eqnarray*}
where the last inequality is by Proposition~\ref{pexp}.
But $|z-w|\le8 e^{-(m-1)\beta}$ implies
\[
m\beta\leq\beta+ \log\frac{8}{|z-w|},
\]
therefore
\[
\log \biggl(\frac{\mathbb{P}[E^z_n \cap E^w_n]}{\mathbb
{P}[E^0_n]^2} \biggr) \le\alpha m \beta \biggl(1+
\frac{\log c}{\alpha\beta} \biggr) \leq \alpha\bigl(1+O(1/\beta)\bigr) \biggl[ \beta+
\log\frac{1}{|z-w|} \biggr].
\]
As $\beta\to\infty$, the error term goes to $0$,
which implies the result.
\end{pf*}

\subsection{Probability of a clockwise loop}
\label{sscw}

In this section, we prove Proposition~\ref{pexp}, lower bounding the
probability that a radial $\mathrm{SLE}_\kappa(\kappa-6)$ process in
$\mathbb{D}$ makes a
clockwise loop within the disk $e^{-\beta}  \mathbb{D}$ before
making any
counterclockwise loop surrounding the origin.

Some notation: for $x\in[0,2\pi]$, we write $\theta^{(x)}_t$ for the
$[0,2\pi]$-valued process of Section~\ref{ssexplorationtree}
started from $\theta^{(x)}_0=x$, evolving according to SDE (\ref
{etheta}) in $(0,2\pi)$ and instantaneously reflecting at the
endpoints $\{0,2\pi\}$. We write $\theta_t\equiv\theta^{(0)}_t$,
and for $a\in[0,2\pi]$ we let $\sigma_a\equiv\inf\{t\dvtx \theta _t=a\}
$, and set $F_t\equiv\{\sigma_{2\pi}>t\}$. For $0<R<1$ and
$\theta_0\in[0,2\pi]$ let $P_R(\theta_0)$ be the probability that a radial
$\mathrm{SLE}_\kappa(\kappa-6)$ in $\mathbb{D}$ with starting configuration
$(w,o)=(1,e^{-i\theta_0})$ and target $0$ makes a clockwise loop inside
the disk $R \mathbb{D}$ surrounding $0$ before making any counterclockwise
loop surrounding $0$. The proposition will be obtained from the
following two lemmas, whose proof we defer.

\begin{lemma}\label{lcond}
There exist $c_0,p_0>0$ such that
\[
\mathbb{P}\bigl[\theta_T \in[c_0,2\pi-c_0]
|F_T\bigr]\ge p_0 \qquad\mbox{for all }T \in[1,\infty).
\]
\end{lemma}

\begin{lemma}\label{lloop}
For any $c_0>0$, we have $\inf_{\theta_0\in[c_0,2\pi
-c_0]}P_R(\theta_0)>0$.
\end{lemma}

\begin{pf*}{Proof of Proposition \ref{pexp}}
Recall that it is natural to parametrize the radial $\mathrm
{SLE}_\kappa(\kappa
-6)$ curve $\eta^0$ targeted at $0$ by capacity: if $U_t$ denotes the
unique connected component of $\mathbb{D}\setminus\eta^0[0,t]$ containing
$0$ and $g_t$ is the unique conformal transformation $U_t\to\mathbb
{D}$ with
$g_t(0)=0$ and $g_t'(0)>0$, then $g_t'(0)=1/\operatorname{CR}(0,U_t)=e^t$.

Assume $\beta\geq\log8$, and
let $\beta'\equiv\beta-\log8\geq0$. Consider the map $g_{\beta
'}\dvtx U_{\beta
'}\to\mathbb{D}$. The left-hand inequality of (\ref{egrowth}) with
$f=g_{\beta'}^{-1}$ gives $|g_{\beta'}(\zeta)|\le4|\zeta| e^{\beta
'}$ for any $\zeta\in U_{\beta'}$.
In particular, $|g_{\beta'}(\zeta)|\le1/2$ for $|\zeta|\le
e^{-\beta'}/8=e^{-\beta}$, so we can apply the right-hand inequality of
(\ref{egrowth}) to find
\[
e^{\beta'} |\zeta| \le 4 \bigl|g_{\beta'}(\zeta)\bigr| \qquad\mbox{when } |\zeta|
\le e^{-\beta}.
\]
Therefore,
the image of $e^{-\beta}  \mathbb{D}$ under $g_{\beta'}$ contains
$R \mathbb{D}$ where
\[
R=\tfrac{1}{4} e^{\beta'-\beta}=\tfrac{1}{32}.
\]
The curve $g_{\beta'}\eta^0$ is distributed as an $\mathrm
{SLE}_\kappa(\kappa-6)$
in $\mathbb{D}$ with starting configuration $(\tilde W_0,\tilde
O_0)=(W_{ \beta
'},W_{ \beta'}  e^{-i\theta_{ \beta'}})$, so for any $c>0$ we have
\[
\mathbb{P}\bigl[E_1^0\bigr] \ge\mathbb{P}[F_{\beta'}]
\mathbb{P}\bigl[\theta_{\beta'}\in[c,2\pi-c] |F_{\beta'}\bigr] \inf_{\theta_0\in[c,2\pi-c]} P_R(\theta_0).
\]
By Proposition~\ref{pssw}, Lemmas~\ref{lcond} and~\ref
{lloop} this expression is $\asymp e^{-\alpha\beta'}$, which gives
the result.
\end{pf*}

The remainder of this subsection is devoted to proving the above lemmas.
We will obtain Lemma~\ref{lcond} as a consequence of the following lemma.

\begin{lemma}
\label{lcondtop}
For any deterministic time $T\geq0$,
\[
\mathbb{P}[\theta_T\le\pi|F_T]\ge1/2.
\]
\end{lemma}

\begin{pf}
Let $S=\pm1$ be a symmetric random sign independent of the process~$\theta_t$, and consider the event $A_T\equiv\{\theta_T<\pi\}\cup
\{\theta_T=\pi, S=1\}$. (The random sign is introduced to handle
the possibility that $\theta_T=\pi$. It follows easily by comparison
with Bessel processes, see, for example, the proof of Lemma~\ref
{lcond}, that $\mathbb{P}[F_T]>0$ and $\mathbb{P}[\theta_T=\pi]=0$
for all
deterministic $T\ge0$, but our proof of Lemma~\ref{lcondtop} can be
applied to any strong Markov continuous process with reflective
symmetry.) By the strong Markov property of $\theta_t$ and the
reflective symmetry across $\pi$ of its drift coefficient,
\[
\mathbb{P}\bigl[A_T^c\bigr] =\mathbb{P}[
\theta_T>\pi]+\tfrac12\mathbb{P}[\theta_T=\pi] =\tfrac12
\mathbb{P}[\sigma_\pi\le T,\theta_T\ne\pi]+\tfrac 12
\mathbb{P}[\theta_T=\pi] \le\tfrac12.
\]
By a similar argument $\mathbb{P}[A_T |F_T^c]\le1/2$. Since $\mathbb
{P}[A_T]\geq
1/2$ is a weighted average of $\mathbb{P}[A_T |F_T^c]\le1/2$
and $\mathbb{P}[A_T |F_T]$, we conclude $\mathbb{P}[A_T |F_T]\ge1/2$,
which implies the lemma.\vadjust{\goodbreak}
\end{pf}

Recall the notation $\theta^{(x)}_t$ defined above. Using the same
driving Brownian motion for any countable collection of processes
$\theta_t^{(x)}$ gives a coupling under which (by continuity and
pathwise uniqueness) the relative order among the processes is
preserved over time.

\begin{pf*}{Proof of Lemma~\ref{lcond}}
For $T\ge1$ and $c_0\in(0,\pi)$, it follows from the Markov property
and Lemma~\ref{lcondtop} that
\begin{eqnarray*}
\mathbb{P} [c_0\leq\theta_T\leq3\pi/2 |F_T
] &\ge &\mathbb{P} [c_0\leq\theta_T\leq3\pi/2,
\theta_{T-1}\le\pi |F_T ]
\\
&\ge&\frac{\mathbb{P}[\theta_{T-1}\le\pi|F_{T-1}]
\inf_{x\le\pi}\mathbb{P} [c_0\leq\theta_1^{(x)}\leq3\pi
/2,F_1^{(x)} ]} {
\mathbb{P}[F_T |F_{T-1}]}
\\
&\ge&\frac12 \inf_{x\le\pi} \mathbb{P} \bigl[\theta_1^{(x)}
\in[c_0,3\pi/2],F_1^{(x)} \bigr]
\\
&\ge&\frac12 \mathbb{P} \Bigl[\theta^{(0)}_1\ge
c_0,\max_{t\le
1}\theta^{(\pi
)}_t
\le3\pi/2 \Bigr],
\end{eqnarray*}
where the last line follows by the coupling described above. Recall
SDE (\ref{etheta}); since $\cot(y/2)/2\sim1/y$ as $y\downarrow0$, by
Girsanov's theorem the process $\theta^{(x)}_t$ before hitting $3\pi
/2$ has law mutually absolutely continuous with respect to that of a
$\sqrt{\kappa} \mathrm{\textsc{bes}}^\delta$ process ($\sqrt {\kappa}$ times a $\delta
$-dimensional Bessel process) started from~$x$, with $\delta\equiv
1+2(\kappa-4)/\kappa$. Note that $\delta>0$ since $\kappa>8/3$. A
$\sqrt{\kappa
} \mathrm{\textsc{bes}}^\delta$ process started from $x\le\pi$
has positive
probability not to hit $3\pi/2$ by time $T$, so $\mathbb{P}[\max_{t\le
1}\theta^{(\pi)}_t\le3\pi/2]>0$. Meanwhile the process $\theta
^{(0)}_t$ before hitting $2\pi$ is mutually absolutely continuous with
respect to another $\sqrt{\kappa} \mathrm{\textsc{bes}}^\delta$
process (started from
zero), so in particular the random variable $\theta^{(0)}_1$ does not
have an atom at $0$ on the event $\{\max_{t\le1}\theta^{(\pi
)}_t\le3\pi/2\}$. Therefore
\[
\lim_{c_0\downarrow0} \mathbb{P} \Bigl[\theta^{(0)}_1
\ge c_0, \max_{t\le1}\theta^{(\pi)}_t
\le3\pi/2 \Bigr] = \mathbb{P} \Bigl[\max_{t\le1}
\theta^{(\pi)}_t\le3\pi/2 \Bigr]>0,
\]
which proves the existence of $c_0,p_0>0$ such that $\mathbb{P}[\theta
_T\in
[c_0,2\pi-c_0]|F_T]\ge p_0$ for all $T\ge1$.
\end{pf*}

\begin{pf*}{Proof of Lemma~\ref{lloop}}
Throughout the proof, let $\eta\equiv\eta_{\theta_0}$ denote a radial
$\mathrm{SLE}_\kappa(\kappa-6)$ process in $\mathbb{D}$ with
starting configuration
$(1,e^{-i\theta_0})$ and target $0$.

We begin by comparing nearby values of $\theta_0$. Let $o=e^{i\theta
_0}$ and
$o'=e^{i\theta_0'}$, where $0<\theta_0,\theta_0'<2\pi$.
The M\"obius transformation
\[
f_{o o'}(\zeta)\equiv \frac{o'}{o} \frac{(o+\bar o'-2)\zeta+(1-o \bar o')} {
(\bar o + o'-2)+(1-\bar o o')\zeta} = \zeta+
\frac{(o'-o) (\zeta-1)^2}{(1-2o+o o')+ (o-o')\zeta}
\]
is the automorphism of $\mathbb{D}$ sending $1$ to $1$, $o$ to $o'$,
and $\bar
o'$ to $\bar o$. Suppose $\theta_0,\theta_0' \in[c_0,2\pi-c_0]$. Then
$|1-2o+o o'| = |2-\bar o-o'|$ is at least $2-\operatorname{Re}(o+o') \ge
2(1-\cos c_0)$, thus bounded away from $0$. From this, it is clear that
if $|\theta_0-\theta_0'|$ is sufficiently small, then the image of
$R
\mathbb{D}$ under $f_{oo'}$ will contain the disk $(R-\epsilon)
\mathbb{D}$. It follows
that $P_R(\theta_0)\ge P_{R-\epsilon}(\theta_0')$ [using the target
invariance of
Proposition~\ref{ptargetinv} since $f_{oo'}(\eta_{\theta_0})$ has
target $f_{oo'}(0)\ne0$]. This reduces the problem of showing
$\inf_{\theta_0\in[c_0,2\pi-c_0]} P_R(\theta_0)>0$ to that of showing
\mbox{$P_R(\theta_0)>0$} for each fixed choice of $0<R<1$ and $\theta
_0\in
[c_0,2\pi-c_0]$.

We prove $P_R(\theta_0)>0$ in two steps which are informally explained in
Figures~\ref{fhook} and \ref{floopclose}.

\begin{figure}

\includegraphics{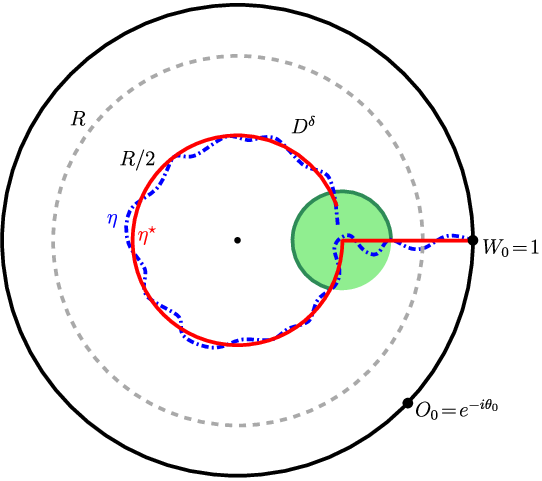}

\caption{Proof of Lemma \protect\ref{lloop}, Step 1: almost clockwise loop. A
radial $\mathrm{SLE}_\kappa(\kappa-6)$
curve $\eta$ (dashed blue line) with starting configuration
$(w,o)=(1,e^{-i\theta_0})$ ($\theta_0\in[c_0,2\pi-c_0]$)
evolves as ordinary chordal $\mathrm{SLE}_\kappa$ from $w$ to $o$, with
(chordal) driving function $W'$ which is $\sqrt{\kappa}$
times a standard Brownian motion. Therefore, $W'$ has positive
probability to be uniformly close to the driving function
$W^\star$ of the hook curve $\eta^\star$. On this event, $\eta$ is
close to $\eta^\star$ in Hausdorff distance and therefore
forms an almost clockwise loop. Write $U$ for the complementary
connected component of the path of $\eta$ so far which contains
$z^\delta$. That $\eta$ closes the clockwise loop with positive
probability is explained in Figure~\protect\ref{floopclose}.}
\label{fhook}
\end{figure}

\textit{Step \textup{1:} Almost clockwise loop.}
The function
\[
f(\zeta) = \frac{i(\zeta-1)(o-1)}{2(\zeta-o)},
\]
conformally maps $\mathbb{D}$ to $\mathbb{H}$ sending $W_0=1$ to $0$
and $O_0=o$ to
$\infty$. Observe
\[
\bigl|f'(\zeta)\bigr| = \biggl\llvert -\frac{i (o-1)^2}{2 (\zeta-o)^2}\biggr\rrvert \ge
\frac{|o^{1/2}-o^{-1/2}|^2}{8} \ge\frac{\sin^2(c_0/2)}{2}\qquad \mbox{for }\zeta\in\mathbb{D},
\]
so the inverse transformation $f^{-1}\dvtx \mathbb{H}\to\mathbb{D}$ is Lipschitz.

Recall from Section~\ref{ssexplorationtree} that up to the stopping
time $\tau_=\equiv\inf\{t\ge0\dvtx W_t=O_t\}$, $\eta$ coincides
(modulo time change) with the exploration tree branch $\eta^{o}$,
which is an ordinary chordal $\mathrm{SLE}_\kappa$ in $\mathbb{D}$
from $W_0=1$ to
$O_0=o$. That is, $f(\eta^o(u))_{u\ge0}$ is a standard chordal
$\mathrm{SLE}
_\kappa$ in $\mathbb{H}$ with associated chordal Loewner driving function
$W'_u=\sqrt{\kappa}B_u$ for $B_u$ a standard Brownian motion, and
$\eta
(t(u))=\eta^o(u)$ for $t(u)\le\tau_=$.

\begin{figure}

\includegraphics{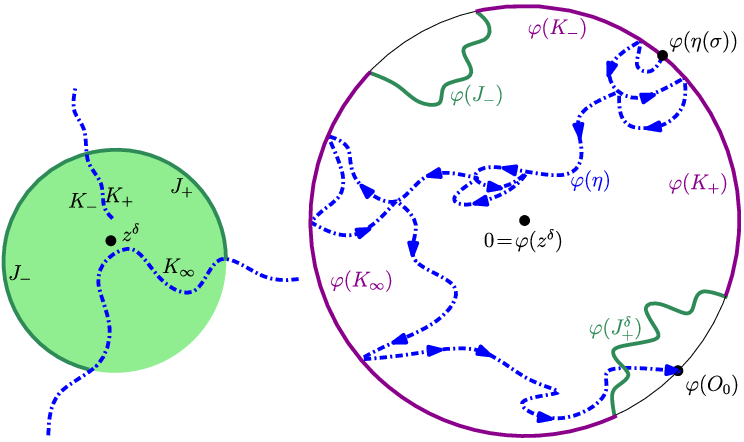}

\caption{Proof of Lemma \protect\ref{lloop}, Step 2: loop closure.
The left panel shows $U$ in a
neighborhood of~$z^\delta$ (see Figure \protect\ref{fhook} for the
notation). The right panel shows the
image of $U$ under the conformal map $\varphi\dvtx U\to\mathbb{D}$ with
$\varphi(z^\delta
)=0$, $\varphi'(z^\delta)>0$.
By conformal invariance of Brownian motion, it follows from
consideration of hitting probabilities of
Brownian motion started from $z^\delta$ in Figure \protect\ref{fhook}
that as $\delta\downarrow0$, $\varphi(J_\pm)$ converge
to points on $\partial\mathbb{D}$ bounded away from one another, and
from the
image of the tip of $\eta$ under $\varphi$.
Loop closure occurs if $\varphi(\eta)$ crosses to the opposing arc
$\varphi
(K_\infty)$ before reaching $\varphi(O_0)$;
this has positive probability for sufficiently small $\delta$ since
$\mathrm{SLE}
_\kappa$ ($4<\kappa<8$) is boundary-intersecting but not boundary-filling.}
\label{floopclose}
\end{figure}

For $0<\delta\ll c_0$, consider the curve $\eta^\star$ in $\mathbb
{D}$ which
travels in a straight line from $1$ to $R/2$, then travels clockwise
along the circle
$\{\zeta\dvtx |\zeta|=R/2\}$ until it reaches $z^\delta\equiv(R/2)
e^{i\delta/2}$. Let $W^\star_u$ be the driving function for $f(\eta
^\star)$ viewed as a chordal Loewner evolution in $\mathbb{H}$,
defined up to
the half-plane capacity $T^\star<\infty$ of $f(\eta^\star)$. By
\cite{MR2129588}, Lemma~4.2, $W^\star_u$ is continuous in $u$, hence
uniformly continuous on $[0,T^\star]$ and thus uniformly approximable
by a piecewise linear function (with finitely many pieces). Since
$W'_u$ is a Brownian motion, for any $\delta'>0$ the event
%
%
\begin{equation}
\label{eBMclose} \sup_{u\le T^\star}\bigl|W'_u-W^\star_u\bigr|
\le\delta'
\end{equation}
occurs with positive probability.
By \cite{MR2129588}, Proposition~4.47, there exists $\delta'>0$ such that if
(\ref{eBMclose}) occurs, then $f(\eta^o[0,T^\star])$ is within
Hausdorff distance $\delta^2\sin^2(c_0/2)/2$ of $f(\eta^\star)$, so
$\eta^o[0,T^\star]$ is within Hausdorff distance $\delta^2$ of $\eta
^\star$. For sufficiently small $\delta$ this implies $t(T^\star
)<\tau_=
$, therefore $\eta[0,t(T^\star)]$ coincides with $\eta^o[0,T^\star
]$. Thus, if we define stopping times
\[
\sigma\equiv\inf\bigl\{t \geq0\dvtx \arg\eta(t)=\delta\bigr\},\qquad \acute\sigma
\equiv\inf\bigl\{t \geq0 \dvtx \operatorname{dist}\bigl(\eta(t),\eta ^\star
\bigr)\ge \delta^2\bigr\},
\]
then we will have $\sigma<\acute\sigma$ on the event
(\ref{eBMclose}).\vadjust{\goodbreak}

\textit{Step \textup{2:} Loop closure.}
On the event $\{\sigma<\acute\sigma\}$, let $\tau$ be the first
time that $\eta$ closes a clockwise loop inside the disk $R  \mathbb{D}$,
and $\acute\tau$ the first time after $\sigma$ that $\eta$
exits $B_{R/4}(R/2)$; the result will follow by showing that
%
%
\begin{equation}
\label{ecloseloop} \liminf_{\delta\downarrow0} \mathbb{P}[\tau<\acute\tau|
\sigma<\acute\sigma]>0.
\end{equation}
Let $U$ denote the unique connected component of $\mathbb{D}\setminus
\eta
[0,\sigma]$ whose closure contains both $0$ and $O_0$. Recall that
$z^\delta\equiv(R/2) e^{i\delta/2}\in U$, and let $\varphi\equiv
\varphi^\delta$
denote the uniformizing map $U\to\mathbb{D}$ with $\varphi(z^\delta
)=0$ and $\varphi
'(z^\delta)>0$. Let $J_{+}$ (resp., $J_{-}$) denote the unique connected
component of $U\cap\partial B_{R/4}(R/2)$ containing the point
$(R/2)+e^{i\pi/4} (R/4)$ (resp., $R/4$): the $J_\pm$ are disjoint
crosscuts\setcounter{footnote}{1}\footnote{A \emph{crosscut} $J$ of a domain $D$ is an open
Jordan arc in $D$ such that $\overline J = J \cup\{a,b\}$ with $a,b\in
\partial D$; $a=b$ is allowed. A crosscut separates the domain into exactly
two components \cite{MR1217706}, Proposition~2.12, and if $\varphi$ is a
conformal map $D\to\mathbb{D}$ then $\varphi J$ is a crosscut of
$\mathbb{D}$ \cite{MR1217706},
Proposition~2.14.} of $U$, and we write $G$ for the
connected component of $U\setminus(J_+ \cup J_-)$
containing $z^\delta$.
The boundary $\partial G$ has a parametrization as a closed curve
$b\dvtx [0,2\pi
]\to\partial G$, oriented counterclockwise with $b(0)=b(2\pi)=\eta
(\sigma
)$.\footnote{The set $A=\partial\mathbb{D}\cup\eta[0,\sigma]\cup
J_\pm$ is
compact, connected, and (since it is a finite union of curves defined
on compact intervals) locally connected. By Torhorst's theorem (see
\cite{MR0507768}, page 285, Problem~1 or \cite{MR0007095}, page 106,
Theorem~2.2), for such $A$, each connected component of $\hat
\mathbb{C}\setminus A$ has a locally connected boundary. In
particular, $\partial
G$ is locally connected, so has a parametrization as a closed curve by
the Hahn--Mazurkiewicz theorem.} We then define times
$0<t_1<t_2<t_3<t_4<2\pi$ such that $b(t_1,t_2)=J_-$ and
$b(t_3,t_4)=J_+$, and write $K_-\equiv b[0,t_1]$, $K_+\equiv b[t_4,2\pi
]$, $K_\infty\equiv b[t_2,t_3]$.

By the conformal Markov property, the probability of $\{\tau <\acute
\tau\}$, conditioned on the path $\eta$ up to time $\sigma
$ on the event $\{\sigma<\acute\sigma\}$, is given by the
probability that a chordal $\mathrm{SLE}_\kappa$ traveling in
$\mathbb{D}$ from $\varphi
(\eta(\sigma))$ to $\varphi(O_0)$ hits $\varphi(K_\infty)$ before
hitting $\varphi(J_\pm)$.\footnote{Here we abuse notation and write
$\varphi S$ for the pre-image of $S$ under the map $\varphi
^{-1}\dvtx \mathbb{D}\to U$
which has a continuous extension to $\overline\mathbb{D}$.} It
follows from
consideration of hitting probabilities of Brownian motion traveling in
$U$ started from $z^\delta$ (using, e.g., the Beurling estimate
\cite{MR2129588}, Theorem~3.76) that as $\delta\downarrow0$, the diameters of the
$\varphi(J_\pm)$ tend to zero while the boundary arcs $\varphi
(K_\infty)$
and $\varphi(K_\pm)$ are all of sizes bounded away from zero. Since
$\mathrm{SLE}_\kappa$ ($4<\kappa<8$) is a.s. boundary-intersecting
but not
boundary-filling (see \cite{MR2129588}, Proposition~6.8) it follows that for
sufficiently small $\delta$ this probability is positive.
\end{pf*}

\subsection{Hausdorff dimension}
\label{sshd}

In this section, we use the second moment estimate Lemma~\ref{lind}
and the lower bound Proposition~\ref{pexp} to deduce the main result
Theorem. \ref{thd}. The argument is standard (see, e.g.,
\cite{MR2642894}, Lemma~3.4) but we give some details here for completeness.

The \emph{$\gamma$-energy} of a Borel measure $\mu$ on a metric
space $(E,d)$ is
\[
I_\gamma(\mu) =\int_{ E}\int_{ E}
\,d(x,y)^{-\gamma} \,d \mu(x) \,d \mu(y).
\]
If there exists a positive Borel measure on $E$ with finite $\gamma
$-energy, then $E$ has Hausdorff dimension bounded below by $\gamma$
(see, e.g., \cite{MR2604525}, Theorem~4.27).

\begin{pf*}{Proof of Theorem~\ref{thd}}
Following the proof of \cite{MR2642894}, Lemma~3.4, we first show that for
any fixed $\varepsilon>0$, there exists with positive probability a nonzero
Borel measure supported on the $\mathrm{CLE}_\kappa$ gasket with finite
$[2-\alpha(1+2\varepsilon)]$-energy, where $\alpha$ is given by
(\ref{alpha}).

For the full range of $\kappa$ we have $\alpha<2$, so we may assume
$\varepsilon$ is sufficiently small that $\alpha(1+\varepsilon)<2$.
Fix $\beta$ large such that $e^\beta/2$ is an integer and
$\varepsilon(\beta
)\le\varepsilon$, with $\varepsilon(\beta)$ as in the statement of
Lemma~\ref
{lind}, and $c\leq e^{\varepsilon\alpha\beta}$, with $c$ as in the
statement of Proposition~\ref{pexp}.

For $z\in\mathbb{C}$ let $S_r(z)\equiv z+[-\frac{r}{2},\frac
{r}{2})\times
[-\frac{r}{2},\frac{r}{2})$ denote the box with side length $r$
centered at $z$, and write $H\equiv S_1(0)\subset\mathbb{D}$. For
$n\ge0$
let $S^z_n\equiv S_{e^{-n\beta}}(z)$. Since $e^\beta/2$ is an integer,
$H$ can be expressed as the disjoint union
\[
H = \bigsqcup_{z\in H_n} S^z_n,\qquad
H_n\equiv \bigl\{e^{-n\beta}(\mathbb{Z}+1/2) \bigr
\}^2 \cap H.
\]
Recall the events $E_n^z$ defined in (\ref{eDzn}).
Define a random measure $\mu_n$ on $H$ by
\[
\mu_n(A) = \int_{ A} \sum
_{z\in H_n} \frac{\mathbf{1}\{E^z_n\}}{\mathbb{P}[E^z_n]} \mathbf{1}\bigl\{z'\in
S^z_n\bigr\} \,d z',\qquad A\subseteq H.
\]
Then $\mathbb{E}[\mu_n(H)]=1$, and
\[
\mathbb{E} \bigl[ \bigl(\mu_n(H) \bigr)^2 \bigr] =
e^{-4n\beta} \sum_{z,w\in H_n} \frac{\mathbb{P}[E^z_n \cap E^w_n]}{\mathbb{P}[E^z_n]\mathbb{P}[E^w_n]}.
\]
The sum over off-diagonal terms is, by Lemma~\ref{lind},
\[
e^{-4n\beta} \mathop{\sum_{z,w\in H_n}}_{z\neq w}
\frac{\mathbb{P}[E^z_n \cap E^w_n]}{\mathbb{P}[E^z_n]\mathbb{P}[E^w_n]} \leq e^{-2n\beta} e^{\alpha\beta(1+\varepsilon)} \mathop{\sum
_{w\in(e^{-n\beta}\mathbb{Z})^2}}_{
|w|<\sqrt{2}} \frac{1}{|w|^{\alpha(1+\varepsilon)}} \preccurlyeq
e^{\alpha\beta(1+\varepsilon)},
\]
using $\alpha(1+\varepsilon)<2$.
By Proposition~\ref{pexp}, the sum over diagonal terms is
\[
e^{-4n\beta} \sum_{z\in H_n} \frac{1}{\mathbb{P}[E^z_n]} \leq
e^{-2n\beta} \bigl(c e^{\alpha\beta}\bigr)^n \leq
e^{-n\beta[2-\alpha(1+\varepsilon)]} \leq1,
\]
therefore $\mathbb{E}[\mu_n(H)^2]\preccurlyeq e^{\alpha\beta
(1+\varepsilon)}$. Similarly,
\begin{eqnarray*}
&&\mathbb{E}\bigl[I_{2-\alpha(1+2\varepsilon)}(\mu_n)\bigr] \\
&&\qquad= \sum
_{z,w\in H_n} \frac{\mathbb{P}[E^z_n\cap E^w_n]}{\mathbb{P}[E^z_n]\mathbb{P}[E^w_n]} \int_{S^z_n}\int
_{S^w_n} \frac{1}{|z'-w'|^{2-\alpha(1+2\varepsilon)}} \,d w' \,d
z'
\\
&&\qquad\preccurlyeq e^{-4n\beta} \sum_{z,w\in H_n}
\frac{\mathbb{P}[E^z_n\cap E^w_n]}{\mathbb{P}[E^z_n]\mathbb{P}[E^w_n]} \biggl\{ e^{n\beta[2-\alpha(1+2\varepsilon)]} \wedge \frac{1} {
\operatorname{dist}(S^z_n,S^w_n)^{2-\alpha(1+2\varepsilon)}} \biggr\}
\\
&&\qquad\preccurlyeq e^{\alpha\beta(1+\varepsilon)} \biggl(e^{-n\beta\alpha\varepsilon} + e^{-4n\beta} \sum
_{z\ne w \in H_n} \frac{1}{|z-w|^{2-\alpha\varepsilon}} \biggr) \preccurlyeq
e^{\alpha\beta(1+\varepsilon)}.
\end{eqnarray*}
The argument of \cite{MR2642894}, Lemma~3.4, then implies that the
$\mathrm{CLE}
_\kappa$ gasket has Hausdorff dimension $\ge2-\alpha(1+2\varepsilon
)$ with
positive probability.

To go from positive probability to probability one, we again make use
of conditional independence in the $\mathrm{CLE}_\kappa$ process.
Recall the
construction of $\mathcal{L}^a_1$, illustrated in Figure~\ref{fcletree} and
described in Section~\ref{sstreesloops}. At the first time $\tau
_{\mathrm{ccw}}
^a$ that $a$ is surrounded by a counterclockwise loop, the loop
$\mathcal{L}
^a_1$ is formed from $\eta^a|_{[\acute\tau_{\mathrm{ccw}}^a,\tau
_{\mathrm{ccw}}^a]}$ together with
an ordinary chordal $\mathrm{SLE}_\kappa$ curve $\widetilde\eta
^a|_{[\acute\tau_{\mathrm{ccw}}
^a,\infty]}$ from $\eta^a(\tau_{\mathrm{ccw}}^a)$ to $\eta
^a(\acute\tau_{\mathrm{ccw}}^a)$ in
the unique connected component $U$ of $\mathbb{D}\setminus\eta
^a[0,\tau_{\mathrm{ccw}}
^a]$ that has both these points on its boundary. Since $\kappa>4$,
this chordal $\mathrm{SLE}_\kappa$ hits both boundary segments
(between the
start and the target) infinitely often, and there are infinitely many
connected components of $U\setminus\widetilde\eta^a[\acute\tau
_{\mathrm{ccw}}^a,\infty]$
which are to the right of the chordal $\mathrm{SLE}_\kappa$. Each connected
component is surrounded by a clockwise loop formed from a segment of
$\widetilde\eta^a$, so by Proposition~\ref{pcle} the components are filled
in by conditionally independent $\mathrm{CLE}_\kappa$ processes.
Since none
of the components is surrounded by a loop of the original $\mathrm
{CLE}_\kappa
$, the gasket of each small $\mathrm{CLE}_\kappa$ is contained within the
gasket of the original $\mathrm{CLE}_\kappa$. Since each of these infinitely
many conditionally independent small gaskets has Hausdorff dimension
$\ge2-\alpha(1+2\varepsilon)$ with positive probability, the
original gasket
has Hausdorff dimension $\ge2-\alpha(1+2\varepsilon)$ almost surely.

Taking $\varepsilon\downarrow0$ implies the theorem.
\end{pf*}

\section*{Acknowledgments}
We thank Steffen Rohde and Scott Sheffield for helpful comments. Part
of the research leading to this article was done at the Mathematical
Sciences Research Institute.


\phantomsection
\bibliographystyle{hmralpha}
%

%



\printaddresses

\end{document}